\documentclass[aima]{preprint}

\usepackage{amssymb}
\usepackage{graphicx}
\usepackage{amsmath}
\usepackage{epsfig}

\newtheorem{theorem}{Theorem}[section]
\newtheorem{corollary}[theorem]{Corollary}
\newtheorem{proposition}[theorem]{Proposition}
\newtheorem{lemma}[theorem]{Lemma}
\newtheorem{remark}[theorem]{Remark}
\newtheorem{definition}[theorem]{Definition}
\newtheorem{example}[theorem]{Example}

\begin{document}

\title{Conjugacy problem for braid groups and Garside groups\footnote{Both authors partially supported by the
European network TMR Sing. Eq. Diff. et Feuill.}}
\author{Nuno Franco\footnote{Partially supported by SFRH/BD/2852/2000.} \\
\begin{tabular}{cc} Dep. de Matem\'atica, CIMA-UE  \hspace{.3cm} & Universit\'e de Bourgogne \\
Universidade de \'Evora  & Laboratoire de Topologie \\
7000-\'Evora (Portugal) & UMR 5584 du CNRS \\
E-mail: nmf@uevora.pt &  B.P. 47870 \\
              &  21078 - Dijon Cedex (France) \\
              & E-mail: nmf@u-bourgogne.fr
\end{tabular}
\and Juan
Gonz\'alez-Meneses\footnote{Partially supported by
BFM-3207.} \\
Dept. de Matem\'atica Aplicada I \\
ETS Arquitectura  \\
Universidad de Sevilla \\
Avda. Reina Mercedes, 2 \\
41012-Sevilla (Spain) \\
E-mail: meneses@us.es}

\maketitle
\sloppy


\begin{subject}
 Primary: 20F36. Secondary: 20F10.
\end{subject}

\vspace{.5cm}
\begin{abstract}
We present a new algorithm to solve the conjugacy problem in Artin
braid groups, which is faster than the one presented by Birman,
Ko and Lee \cite{B-K-L}. This algorithm can be applied not only to
braid groups, but to all {\em Garside groups} (which include
finite type Artin groups and torus knot groups among others).
\end{abstract}

\section{Introduction}


Given a group $G$, the conjugacy problem in $G$ consists on
finding an algorithm which, given $a,b\in G,$ determines if there
exists $c\in G$ such that $a=c^{-1}bc.$ Sometimes one also needs
to compute $c,$ for instance, when one tries to attack
cryptosystems based on conjugacy in $G$ (\cite {A-A-G},
\cite{K-L}).

We are mainly interested in Artin braid groups, which are defined, for $%
n\geq 2$, by the following presentation:
\begin{equation}
B_{n}=\left\langle
\sigma _{1},\sigma _{2},\ldots ,\sigma _{n-1}  \left|
\begin{array}{ll}
\sigma _{i}\sigma _{j}=\sigma _{j}\sigma _{i} & (|i-j|\geq 2) \\
\sigma _{i}\sigma _{i+1}\sigma _{i}=\sigma _{i+1}\sigma _{i}\sigma _{i+1} &
(1 \leq i \leq n-2)
\end{array}
\right.
\right\rangle   \label{presen}
\end{equation}

 The first conjugacy algorithm for braid groups was given by Garside
\cite{G}. It was improved by Elrifai and Morton \cite{M} and, more
recently, by Birman, Ko and Lee (\cite{B-K-L} and \cite{BKL2}).

In all these algorithms, one of the key points is the existence of a
finite set $S\subset B_{n}$, whose elements are called
\textit{simple elements}, verifying some suitable properties
(we will be more precise later). One of the main disadvantages
is the size of $S$, which is always greater than $3^n$.

In this paper we will show how one can avoid this problem by defining
some small subsets of $S$, whose size is smaller than $n-1$.
Their elements will be called {\em minimal simple elements}.
Unlike $S$, these sets of minimal simple elements are not unique for every group:
The suitable set of minimal simple elements must be recomputed many times in our
algorithm. Nevertheless, we will see that it is much faster to compute and use these
very small subsets, than to use the whole $S$ all the time.

For instance, the known upper bound for the complexity of the Birman-Ko-Lee algorithm, to decide
wether two braids $a$ and $b$ are conjugated in $B_n$, is
$O( kl^{2}n3^{n}) $ (where $k$ is a number that will be explained later, and $l$ is the
maximum of the word lengths of $a$ and $b$).  An upper bound for the complexity of our
algorithm for $B_n$ is $O( kl^2n^4)$.

Let us mention that our algorithm, as well as the previous ones, also computes the
element $c\in B_n$ such that $a=c^{-1}bc$.
Moreover, since our construction relies on the existence of simple
elements and their basic properties, we can extend our results to
a much larger class of groups, called \textit{Garside groups}.
They were introduced by Dehornoy and Paris \cite{DP}.
At the origin, these groups were called {\it small Gaussian groups}, but
there has been a convention to call them {\it Garside groups}.
They include, besides Artin braid groups, spherical (finite type)
Artin groups, torus knot groups and others.

 One final remark: one important property of Garside groups
is the existence of embedable monoids (for instance the monoid of
positive braids, $B_n^+$, which embeds in $B_n$). The conjugacy
class of an element $a$ in such a monoid is known to be a finite
set, $C^+(a)$. We will also show how to compute $C^+(a)$, using
the techniques mentioned above.

 This paper is structured as follows:
In Section~2, we give a brief introduction to Garside monoids and groups;
In Section~3, the known algorithms mentioned in this introduction are detailed;
We introduce the minimal simple elements in Section~4, and
in Section~5 we present our algorithms in detail;
Complexity issues are treated in Section~6 and, finally,
some effective computations are described in Section~7.

\section{Garside monoids and groups}

The results contained in this section are well known, and can be found in
\cite{G}, \cite{M}, \cite{T}, \cite{B-K-L}, \cite{DP}, \cite{D} and \cite{Pi}.
We will define the {\em Garside monoids} and {\em Garside groups}, and explain
some basic properties.

 Given a cancellative monoid $M$, with no invertible elements, we can define two different partial
orders on its elements, $\prec $ and $\succ$. Given $a,b\in M$, we say that $a\prec b$ $\left(
b\succ a\right) $ if there exists $c\in M$ such that
$ac=b$ $\left( b=ca\right) $, and we say that $a$ is a left
(right) divisor of $b$.

In this situation, we can naturally define the (left or right) {\em least common
multiple} and {\em greatest common divisor} of two elements.
Given $a,b\in M$, we denote by $a\vee b$ the left lcm of $a$ and $b$, if it
exists. That is, a minimal element (with respect to $\prec$) such that
$a\prec a\vee b$ and $b\prec a\vee b.$ We denote by $a\wedge b$ the left gcd
of $a$ and $b$, if it exists. That is, a maximal element (with respect to $\prec
),$ such that $a\wedge b\prec a$ and $a\wedge b\prec b$.

\begin{definition}\label{defatom}
 Let $M$ be a monoid. We say that $x\in M$ is an {\em atom} if $x\neq 1$ and if
$x=yz$ implies $y=1$ or $z=1$. $M$ is said to be an {\em atomic monoid} if it is generated
by its atoms and, moreover, for every $x\in M$, there exists an integer $N_x>0$ such that $x$
cannot be written as a product of more than $N_x$ atoms.
\end{definition}

\begin{definition}
 We say that a monoid $M$ is a {\em Gaussian monoid} if it is atomic, (left and right)
cancellative, and if every pair of elements in $M$ admits a (left and right) lcm and
a (left and right) gcd
\end{definition}

\begin{definition}
 A {\em Garside monoid} is a Gaussian monoid which has a {\em Garside element}. A {\em Garside
element} is an element $\Delta\in M$ whose left divisors coincide with their right divisors,
they form a finite set, and they generate $M$.
\end{definition}

\begin{definition}
 The left (and right) divisors of $\Delta$ in a Garside monoid $M$ are called
\textit{simple elements}. The (finite) set of simple elements is denoted by $S$.
\end{definition}

It is known that every Garside monoid admits a group of fractions. So we have:

\begin{definition}
 A group $G$ is called a {\em Garside group} if it is the group of fractions of a Garside monoid.
\end{definition}

The main example of a Garside monoid (actually the monoid studied by Garside) is the Artin braid
monoid on $n$ strands, $B_{n}^{+}$. It is defined by
Presentation~(\ref{presen}), considered as a presentation for a monoid.
Its group of fractions is the braid group $B_{n}$, and Garside \cite{G} showed that $B_{n}^{+}\subset B_{n}.$
Actually, every Garside monoid embeds into its corresponding Garside group~\cite{DP}.

\begin{figure}[ht]
\centerline{\includegraphics{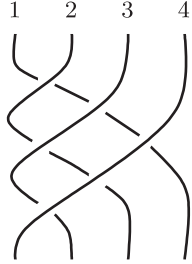}}
\caption{The Garside element $\Delta\in B_4^+$.}\label{Delta4}
\end{figure}

The classical choice of a Garside element for $B_n^+$ is the following: $\Delta =\left( \sigma _{1}\sigma
_{2}\cdots \sigma _{n-1}\right) \left( \sigma _{1}\sigma
_{2}\cdots \sigma _{n-2}\right) \cdots \left( \sigma _{1}\sigma
_{2}\right) \sigma _{1}.$ It can be defined as the positive braid (braid in $B_{n}^{+}$) in which any two
strands cross {\em exactly} once (where, as usual, $\sigma_i$
represents a crossing of the strands in positions $i$ and
$i+1$). It is represented in Figure~\ref{Delta4} for $n=4$.
The simple elements in this case are the positive braids in which any two strands cross {\em at most} once.
Then one has $\#(S)=n!$

Another important example of Garside monoid is the Birman-Ko-Lee
monoid \cite {B-K-L}, which has the following presentation:
\begin{equation}
BKL_{n}^{+} =  \left\langle
  a_{ts} (n\geq t>s\geq 1)
 \left|
\begin{array}{l}
 a_{ts}a_{rq}=a_{rq}a_{ts} \text{ if }\left( t-r\right) \left( t-q\right) \left( s-r\right) \left(
s-q\right) >0 \\
a_{ts}a_{sr}=a_{tr}a_{ts}=a_{sr}a_{tr} \;\text{ where } n\geq t>s>r\geq 1
\end{array}
\right . \right\rangle
\label{presen2}
\end{equation}


Its group of fractions is again  the braid group $B_{n}$. The usual Garside element in
$BKL_{n}^{+}$ is $\delta =a_{n,n-1}a_{n-1,n-2}\cdots a_{2,1}.$ The
advantage of this monoid with respect to $B_n^+$ is that
$\#(S)=\mathcal{C}_{n}$, where
$\mathcal{C}_{n}=\frac{\left( 2n\right) !}{n!\left(
n+1\right) !}<4^{n}$ is the $n^{th}$ Catalan number. Hence, the number of simple elements
is much smaller in this case, but it is still quite big, since $\mathcal{C}_{n}>3^n$.
Notice also that $\left| \delta \right| =n-1,$ while in
$B_n^+$, $\left| \Delta \right| =\frac{n\left( n-1\right) }{2}.$

As we mentioned before, there are other examples of Garside groups, such as finite type Artin groups,
or torus knot groups (see~\cite{Pi} to find more examples of Garside groups).

From now on, $M$ will
denote a Garside monoid, $G$ its group of fractions and $\Delta $ the corresponding Garside
element. Since $M\subset G$, we will refer to the elements in $M$
as the {\em positive} elements of $G$.

\begin{figure}[ht]
 \centerline{\includegraphics{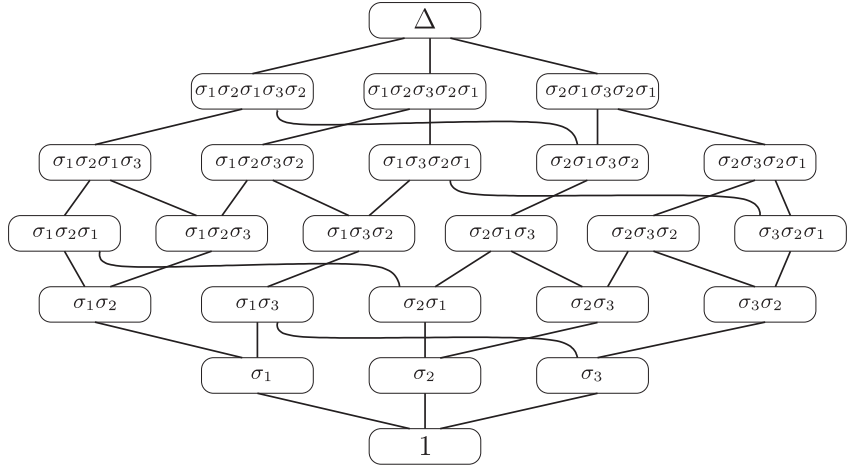}}
\caption{The lattice of simple elements in
$B_4^+$.}\label{lattice}
\end{figure}

From the existence of left lcm's and gcd's, it follows that
$\left( M,\prec \right) $ has a lattice structure, and
$S$ becomes a finite sublattice with minimum $1$ and maximum
$\Delta .$ See in Figure~\ref{lattice} the Hasse diagram of the lattice of $S$ in
$B_{4}^{+},$ where the lines represent left divisibility (from
bottom to top). The analogous properties are also verified by $\succ $.

\begin{definition}
For $a\in M$ we define $LM\left( a\right) \in S$ as the maximal
simple left divisor of $a,$ that is, $LM\left( a\right) =\Delta
\wedge a.$ We also define $RM\left( a\right) $ as the maximal
simple right divisor of $a$.
\end{definition}

\begin{proposition} [\cite{G}]
For $a\in G$, there exists a unique decomposition $a=\Delta
^{p}a_{1}\cdots a_{l}$, called {\em left normal form} of $a$, where:

\begin{enumerate}
\item  $p=\max \left\{ r\in \mathbb{Z}:\Delta ^{-r}a\in M\right\} $
 \quad (hence $a_1\cdots a_l\in M$).

\item  $a_{i}=LM\left( a_{i}\cdots a_{l}\right) \in
S\backslash \{\Delta,1\}$, for all $i=1,...,l.$
\end{enumerate}
Symmetrically, one defines the {\em right normal form} of $a\in G,$
using $RM$.
\end{proposition}

Sometimes, if we are dealing with elements in $M$ and it does not lead to confusion,
we will say that an element $w=w_{1}\cdots w_{t}\in M$ is in left
normal form to express that $w_i\in S\backslash \{1\}$ for all $i$ and, for some $p\geq 0$,
the normal form of $w$ is $\Delta^p w_{p+1}\cdots w_t$.

Later we will use these technical results:

\begin{lemma} [\cite{Mi}, Prop. 2.1]
\label{corLM}
  Let $w_{1}\cdots w_{t}\in M$ be in left
normal form, and $x_1\cdots x_t\in M$ in right normal form. For
every $v\in M$, one has $LM\left( vw_{1}\cdots w_{t}\right)
=LM\left( vw_{1}\right)$ and $RM \left( x_{1}\cdots x_{t}v\right)
=RM\left( x_{t}v\right)$.
\end{lemma}

\begin{lemma} [\cite{Mi}, Prop. 5.3]
\label{lemRNF}
Let $w=w_1\cdots w_t\in M$ be written in right normal form. If we write $w$ in any other way
as a product of $t$ simple elements, $w=u_1\cdots u_t$, then $w_1\prec u_1$.
\end{lemma}

\begin{lemma}[\cite{C}, 3.1] \label{lemcharney}
Let $w=w_1 \cdots w_t\in M$ be written in right normal form, and
let $s\in S$. Then we can decompose $w_i=w_i'w_i''$, for all $i$, in such a way that the right
normal form of $ws$ is $(w_1')(w_1''w_2')\cdots (w_{t-1}''w_t')(w_t''s)$ if it has $t+1$ factors,
or $(w_1 w_2')\cdots (w_{t-1}''w_t')(w_t''s)$ if it has $t$ factors.
\end{lemma}

\begin{corollary}
\label{mainSSS}
Let $w=w_1\cdots w_t\in M$ be written in right normal form. Let $s\in S$ and suppose that we can write
$w s$ as a product of $t$ simple elements, that is,
$w_{1}\cdots w_{t}s=u_1\cdots u_t$. Then $w_{1}\prec u_{1}$.
\end{corollary}

\begin{proof}
Since $w s$ can be written as a product of $t$ simple elements, then its right normal form has
$t$ factors, say $v_1\cdots v_t$. By Lemma~\ref{lemcharney},
$w_1\prec v_1$, and by Lemma~\ref{lemRNF} $v_1\prec u_1$, so the result follows.
\end{proof}

We end this section with a last property of Garside groups: There is a power
of their Garside element which belongs to the center. For instance, in $B_{n}$ the element
$\Delta^{2}=\delta ^{n}$ generates the center of $B_n$.

\section{Known algorithms for the conjugacy problem.}

We present here the Elrifai-Morton algorithm for the conjugacy problem in braid groups \cite{M},
which is also valid for Garside groups, as can be seen in \cite{Pi2}.

It goes as follows: for every element $a\in G$, it computes a
finite subset $C^{sum}(a)$ of the conjugacy class of $a$. This set
is shown to be independent of $a$, so it is an invariant of its
conjugacy class. Therefore, two elements $a$ and $b$ are
conjugated if and only if $C^{sum}(a)=C^{sum}(b)$.

Let us explain the algorithm in more detail.

\subsection{Definition of $\mathbf{C^{\geq m}(a)}$ and $\mathbf{C^{sum}(a)}$}

\begin{proposition}{\em \cite{M,Pi2}}
 Let $a=\Delta ^{p}a_{1}\cdots a_{l}\in G$ be in left normal form.
Then the right normal form of $a$ is as follows:
$a=x_1\cdots x_l \Delta^p$, where $l$ and $p$ are the same as above.
\end{proposition}

\begin{definition}
Let $a=\Delta ^{p}a_{1}\cdots a_{l}\in G$ be in left normal form.
We define the {\em infimum}, {\em supremum} and {\em canonical
length} of $a$, respectively, by $\inf \left( a\right) =p,$ $\sup
\left( a\right) =p+l$, and $\left\| a\right\| =l$ .
\end{definition}

\begin{definition}
Let $a\in G$ and denote by $C\left( a\right) $ the conjugacy class
of $a$. We define the {\em summit infimum}, the {\em summit
supremum} and the {\em summit length} of $a$ as, respectively,
$\max \left\{ \inf \left( x\right) :x\in C\left( a\right) \right\}
$, $\min \left\{ \sup \left( x\right) :x\in C\left( a\right)
\right\} $ and $\min \left\{ \left\| x\right\| :x\in C\left(
a\right) \right\} .$
\end{definition}

\begin{definition}
Let $a \in G$.
\begin{enumerate}
\item For every integer $m$, we define $C^{\geq m}(a)=\left\{v\in C(a): \;
\inf(v)\geq m\right\}$.

\item We define the {\em summit class} of $a$, $C^{sum}\left( a\right)$, as the subset
of $C(a)$ containing all elements of minimal canonical length.
\end{enumerate}
\end{definition}

\noindent {\bf Remarks:}

1. One has $C^{\geq 0}(a)=C(a)\cap M=C^+(a)$.

2. In~\cite{M}, $C^{sum}\left( a\right)$ is called
the {\em Super Summit Set}.

\begin{proposition}{\em \cite{M,Pi2}}\label{prosum}
For every $b\in C^{sum}(a)$, the infimum, supremum and canonical
length of $b$ are equal, respectively, to the summit infimum, the
summit supremum and the summit length of $a$.
\end{proposition}

 It is known that $C^{\geq m}(a)$ and $C^{sum}\left( a\right)$ are finite
sets. Moreover, by Proposition~\ref{prosum}, if
$C^{\geq m}\left( a\right)\neq \phi$, then $C^{sum}\left(
a\right)\subset C^{\geq m}\left( a\right)$.

\subsection{Cycling and decycling}

 Let $\tau :G\rightarrow G$ be the automorphism defined
by $\tau(a)= \Delta^{-1} a \Delta$. The restriction of $\tau$ to $S$ is a
bijection $\tau: S\rightarrow S$.

\begin{definition}
Let $a=\Delta ^{p}a_{1}\cdots a_{l}\in G$ be written in left
normal form. The functions {\em cycling} and {\em decycling} are the
maps $\mathbf{c}$ and $\mathbf{d},$ from $G$ to itself,
defined by:
\begin{eqnarray*}
\mathbf{c}\left( a\right)  &=&\Delta ^{p}a_{2}\cdots a_{l}\tau ^{-p}\left(
a_{1}\right) ; \\
\mathbf{d}\left( a\right)  &=&\Delta ^{p}\tau ^{p}\left( a_{l}\right)
a_{1}\cdots a_{l-1}.
\end{eqnarray*}
\end{definition}

Notice that $\mathbf{c}\left( a\right) $ and $\mathbf{d}\left(
a\right) $ are conjugates of $a$. Furthermore, for every $a\in G$,
$\inf(a)\leq \inf(\mathbf{c}(a))$ and
$\sup(a)\geq \sup(\mathbf{d}(a))$.

Suppose that we have an element $a\in G$, such that $\inf(a)$ is not equal
to the summit infimum of $a$. Then we can try to increase the infimum by repeated
cycling. By~\cite{M} (and~\cite{Pi2}), this always works: there exists a positive
integer $k$ such that $\inf(\mathbf{c}^k(a))>\inf(a)$. We know a bound for this
integer $k$ only for some special Garside monoids and groups: If $M$ is {\em homogeneous},
i.e. it has only homogeneous relations (for instance, if $M$ is $B_n^+$ or $BKL_n^+$),
then every two words representing an element $a\in M$ have the same length, denoted $|a|$.
It is shown in~\cite{BKL2} that, in this case, $k<|\Delta|$.

Therefore, by repeated cycling, we can conjugate $a$ to another element $\widehat{a}$ of maximal
infimum. Even if $M$ is not homogeneous, we know that we reached the summit infimum when we enter into a
loop: at some point $\mathbf{c}^k(v)=v$ for some $v$ conjugated to $a$. This always happens since the
set $C^{\geq m}(a)$ is finite for every $m$, in particular for the summit infimum.

Once $\widehat{a}$ is obtained, we can try to decrease its supremum by repeated decycling.
By~\cite{M} (and~\cite{Pi2}), this also works: either we enter into a loop, and then the supremum is
minimal, or there exists an integer $k$ such that $\sup(\mathbf{d}^k(\widehat{a}))
<\sup(\widehat{a})$. Again
by~\cite{BKL2}, $k<|\Delta|$ in homogeneous monoids.

Therefore, using repeated cycling and decycling a finite number of times, one
obtains an element $\widetilde{a} \in C^{sum}\left( a\right)$.
And, if $M$ is homogeneous, this can be done in polynomial time in $|a|$.

\subsection{The Elrifai-Morton algorithm}

Once that we obtained an element $\widetilde{a} \in C^{sum}(a)$,
we can construct the whole $C^{sum}\left( a\right)$, by using the next result:

\begin{proposition}
{\em \cite{M,Pi2}}\label{prop17} For $u,v$ conjugate
elements in $C^{sum}(a)$ (resp. $C^{\geq m}( a)$), there
exists a sequence $u=$ $u_{1},u_{2},...,u_{k}=v$ of elements in
$C^{sum}(a)$ (resp. $C^{\geq m}( a)$) such that, for
$i=1,\ldots,k-1$, $u_{i} $ and $\ u_{i+1}$ are conjugated by an
element in $S$.
\end{proposition}

 The Elrifai-Morton algorithm does the following: Given $a,b\in
G$ it computes, using cyclings and decyclings, $\widetilde{a}\in C^{sum}(a)$ and
$\widetilde{b}\in C^{sum}(b)$. Then it defines $V_1=\{\widetilde{a}\}$
and it computes, by recurrence,
$$
V_i=\{s^{-1}vs;\; s\in S, v\in V_{i-1}\}\cap C^{sum}(a).
$$
Since $1\in S$, this creates an ascending chain of subsets of $C^{sum}(a)$. By
the above proposition, one has $V_k=V_{k+1}$ for some $k$, and then
$V_k=C^{sum}(a)$. Hence, when the chain stabilises, the whole $C^{sum}(a)$ has been computed.
Then $a$ and $b$ are conjugated if and only if $\widetilde{b}\in C^{sum}\left(a\right)$.

\begin{remark}
This algorithm can be modified to compute $C^{\geq m}(a)$ for $a\in M$ and $m\in \mathbb{Z}$.
We just need to replace $C^{sum}(a)$ by $C^{\geq m}(a)$ in the above discussion.
\end{remark}

Notice that $C^{sum}(a)$ (resp. $C^{\geq m}(a)$) is computed at the cost of conjugating
every element in $C^{sum}(a)$ (resp. $C^{\geq m}(a)$) by every element in $S$. All these sets are
quite big, and this makes the algorithm to be slow. In what follows, we will
get rid of the problem caused by the size of $S$, using the {\em minimal simple elements}.

\section{Minimal simple elements}\label{secmin}

 In this section we shall define some very small subsets of $S$,
which will enable us to compute $C^{\geq m}(a)$ and $C^{sum}(a)$, for $a\in G$, much
faster than the previous algorithms.

Recall the definition of the partial order $\prec$ in $M$.

\begin{definition}
Let ${\cal P}$ be a property for simple elements.
We denote by $S_{\cal P}$ the set of simple elements satisfying ${\cal P}$.
The {\em set of minimal simple elements} for ${\cal P}$,
$\min (S_{\cal P})$, is the set of minimal elements (with respect to $\prec$) in
$S_{\cal P}$.
\end{definition}

 We shall enforce ${\cal P}$ to be {\em closed under g.c.d}, that is, if $s_1,s_2\in S_{\cal P}$
then $s_1\wedge s_2 \in S_{\cal P}$. Let us see that, under this assumption, the set
$\min(S_{\cal P})$ turns to be very small.
For every atom $x\in M$, let $mult(x)=\{ s\in S; \; x\prec s\}$.

\begin{lemma}\label{lemrhox}
 Suppose that ${\cal P}$ is closed under gcd, and let $x$ be an atom of $M$.
If the set $S_{\cal P} \cap mult(x)$ is non-empty, then it has a unique
minimal element, that we denote $\rho_x$.
\end{lemma}

\begin{proof}
 Suppose that there are two distinct minimal elements $s_1,s_2\in S_{\cal P} \cap mult(x)$.
Since $s_1,s_2\in S_{\cal P}$, then  $s_1\wedge s_2\in S_{\cal P}$. Moreover, since $x$ divides
$s_1$ and $s_2$, it also divides $s_1\wedge s_2$. Therefore $s_1\wedge s_2\in S_{\cal P} \cap mult(x)$,
so $s_1$ and $s_2$ cannot be both minimal.
\end{proof}

\begin{corollary}\label{coratoms}
 Suppose that $M$ has $m$ atoms. If ${\cal P}$ is
closed under gcd, then $\#(\min(S_{\cal P}))\leq m$.
\end{corollary}

\begin{proof}
 Notice that every element in $M$ must be divisible by an atom. Take $s\in \min(S_{\cal P})$
and consider an atom $x\prec s$. Since $s$ is minimal in $S_{\cal P}$, it is also minimal in
$S_{\cal P} \cap mult(x)$. Hence $s=\rho_x$. Therefore
$$
\min(S_{\cal P})\subset \{\rho_x: \; x \mbox{ is an atom}\}
$$
and the result follows.
\end{proof}

\begin{example}
 In $B_n^+$ there are $n-1$ atoms, namely $\sigma_1,\ldots, \sigma_{n-1}$. Therefore, if ${\cal P}$
is a property closed under gcd, then $\min(S_{\cal P})$ has at most $n-1$ elements, while $\#(S)=n!$
\end{example}

\begin{example}
 In $BKL_n^+$ there are $\frac{n(n-1)}{2}$ atoms (the generators in Presentation~\ref{presen2}). Hence,
if ${\cal P}$ is a property closed under gcd, then $\#(\min(S_{\cal P}))\leq \frac{n(n-1)}{2}$,
while $\#(S)={\cal C}_n>3^n$.
\end{example}

We must now define some suitable properties, closed under gcd, that will allow us to compute
$C^{\geq m}(a)$ and $C^{sum}(a)$, for $a\in G$. These properties will depend on some given elements
in $M$, so we will have an infinite number of properties, each one corresponding to a set
of minimal simple elements.

\subsection{Minimal simple elements to compute $\mathbf{C^{\geq m}(a)}$}

\begin{definition}
 Let $a\in G$ and $v\in C^{\geq m}(a)$, for some $m\in \mathbb{Z}$. We will say that a simple
element $s$ satisfies the property ${\cal P}_v^{\geq m}$ if it conjugates $v$ to an element in
$C^{\geq m}(a)$, that is, $s^{-1}vs\in C^{\geq m}(a)$.
\end{definition}

\begin{proposition} \label{propcarac}(Caracterization of elements satisfying ${\cal P}_v^{\geq m}$).
If $v\in C^{\geq m}(a)$, one can write $v=\Delta^{m}w$, where $w\in M$. Then a simple element
$s$ satisfies the property ${\cal P}_v^{\geq m}$ if and only if $\tau^m(s)\prec ws$.
\end{proposition}

\begin{proof}
 The first assertion comes from the definition of infimum. Let then $v=\Delta^m w$, where $w\in M$,
and let $s\in S$. One has $s^{-1} v s = s^{-1} \Delta^m w s = \Delta^m \tau^m(s^{-1}) w s=
\Delta^m (\tau^m(s))^{-1} w s$. Hence, $s$ satisfies ${\cal P}_v^{\geq m}$ if and only if
$(\tau^m(s))^{-1} w s\in M$,
that is, $\tau^m(s)\prec ws$.
\end{proof}

\begin{proposition}\label{prop+gcd}
 For every $v\in M$ and every $m\in \mathbb{Z}$, the property ${\cal P}_v^{\geq m}$ is closed under gcd.
\end{proposition}

\begin{proof}
 Suppose that $s_1$ and $s_2$ satisfy ${\cal P}_v^{\geq m}$, and let $s=s_1\wedge s_2$. Notice that
 $\tau$ preserves gcd's, since it preserves left divisibility. Hence $\tau(s)=\tau(s_1)\wedge \tau(s_2)$,
 and thus $\tau^m(s)=\tau^m(s_1)\wedge \tau^m(s_2)$.

One has $\tau^m(s)\prec \tau^m(s_1)\prec vs_1$  and $\tau^m(s) \prec \tau^m(s_2)\prec vs_2$. But it is easy
to show that, for every $v\in M$, $vs_1\wedge vs_2=vs$. Hence, since $\tau^m(s)$ divides $vs_1$ and $vs_2$
then it divides its gcd, i.e. $\tau^m(s)\prec vs$. Therefore, $s$ satisfies ${\cal P}_v^{\geq m}$, and the
result follows.
\end{proof}

\begin{definition}
 For every $v\in C^{\geq m}(a)$, we define $S_v^{\geq m}=min(S_{{\cal P}_v^{\geq m}})$.
That is, $S_v^{\geq m}$ is the set of minimal simple elements (with respect to $\prec$)
among those who conjugate $v$ to an element in $C^{\geq m}(a)$.
\end{definition}

Notice that, by Corollary~\ref{coratoms} and Proposition~\ref{prop+gcd}, the cardinal of $S_v^{\geq m}$
for every $v\in C^{\geq m}(a)$ is no bigger than the number of atoms in $M$. Moreover, we have the following
result, analogous to Proposition~\ref{prop17}.

\begin{proposition}
\label{cormainA1}Given $u,v\in C^{\geq m}(a)$ for some $a\in G$, there exists a sequence $%
u=u_{1}, u_2,...,u_{k}=v$ of elements in $C^{\geq m}(a)$ such that, for
$i=1,...,k-1$, the elements $u_i$ and $u_{i+1}$ are conjugated by an element in $S_{u_i}^{\geq m}$.
\end{proposition}

\begin{proof}
Just notice that any left or right divisor of a simple element is
also a simple element, and then decompose every simple element in
the sequence given by Proposition~\ref{prop17} into a product of
minimal ones.
\end{proof}

This result implies that, in order to compute $C^{\geq m}(a)$ for $a\in M$, it
suffices to conjugate every $v\in C^{\geq m}(a)$ by the elements
in the small set $S_{v}^{\geq m}$.

\subsection{Minimal simple elements to compute $\mathbf{C^{sum}(a)}$}

\begin{definition}
 Let $a\in G$, and let $v\in C^{sum}(a)$. We will say that a simple element $s$ satisfies the
property ${\cal P}_v^{sum}$ if it conjugates $v$ to an element in $C^{sum}(a)$.
In other words, if the canonical length of $s^{-1}vs$ is equal to the canonical length of $v$
(which is the summit length of $a$).
\end{definition}

\begin{proposition}\label{corSSS}
 For every $v\in C^{sum}(a)$, the property ${\cal P}_v^{sum}$ is closed under gcd.
\end{proposition}

\begin{proof}
Let $s_1$ and $s_2$ be two simple elements satisfying ${\cal P}_v^{sum}$, and
denote $s=s_1\wedge s_2$. Write $s_{i}=s r_{i}$ for $i=1,2$, thus $r_1\wedge r_2=1$.

Suppose that $\inf(v)=p$ and $\left\| v \right\|=t$. Then $v=\Delta^p v'$, where
$v'\in M$ and we can write $v'$ as a product of $t$ simple elements (but not less).
Since $s_1$ satisfies ${\cal P}_v^{sum}$, one has
$ s_1^{-1} v s_1 = s_1^{-1}\Delta^p \:v' s_1 = \Delta^p\: \tau^p(s_1^{-1}) \:v' s_1 =
  \Delta^p \:(\tau^p(s_1))^{-1} v' s_1,$
where $(\tau^p(s_1))^{-1} v' s_1\in M$ and we can write it as a product of $t$ simple elements,
say $x_1\cdots x_t$. The same happens for $(\tau^p(s_2))^{-1} v' s_2\in M$.

 Now consider $s^{-1}vs$. By Proposition~\ref{prop+gcd} it belongs to $C^{\geq p}(a)$, that is,
$(\tau^p(s))^{-1} v' s\in M$. We must show that we can write this element as a product of $t$ simple
elements. Suppose this is not true, and write $(\tau^p(s))^{-1} v' s=z_1\cdots z_{t+1}$
in right normal form (it has no more than $t+1$ factors since it is a right divisor of $v's$ which has
$t+1$ factors). One has
$ x_1\cdots x_t= (\tau^p(s_1))^{-1} v' s_1 = (\tau^p(r_1))^{-1} (\tau^p(s))^{-1} v' s r_1=
(\tau^p(r_1))^{-1} z_1\cdots z_{t+1} r_1.$
Hence, $z_1\cdots z_{t+1} r_1=\tau^p(r_1) x_1\cdots x_t$, and $z_1\cdots z_{t+1}$ is in right normal form.
Then by Corollary~\ref{mainSSS}, $z_1\prec \tau^p(r_1)$. In the same way, $z_1\prec \tau^p(r_2)$.
Therefore $z_1\prec \tau^p(r_1)\wedge \tau^p(r_2)= \tau^p(r_1\wedge r_2)= \tau^p(1)=1$.
A contradiction.
\end{proof}

\begin{definition}
 For every $v\in C^{sum}(a)$, we define $S_v^{sum}=min(S_{{\cal P}_v^{sum}})$. That is, $S_v^{sum}$
is the set of minimal simple elements (with respect to $\prec$) among those who conjugate
$v$ to an element in $C^{sum}(a)$.
\end{definition}

As before, by Corollary~\ref{coratoms} and Proposition~\ref{corSSS}, the cardinal of $S_v^{sum}$
for every $v\in M$ is no bigger than the number of atoms in $M$. Furthermore, we can adjust the algorithm
by Elrifai-Morton to these new sets, since we have the following
result, analogous to Propositions~\ref{prop17} and \ref{cormainA1}.

\begin{proposition}
For $u,v$ conjugate elements in $C^{sum}\left( a\right) $, there
exists a sequence $u=u_{1},...,u_{k}=v$ of elements in
$C^{sum}\left( a\right) $ such that, for $i=1,...,k-1$, the
elements $u_i$ and $u_{i+1}$ are conjugated by an element in
$S_{u_{i}}^{sum}$.
\end{proposition}

 The proof of this result parallels that of Proposition~\ref{cormainA1}.
It implies that, in order to compute $C^{sum}\left( a\right) $ for $a\in G$, it
suffices to conjugate every $v\in C^{sum}\left( a\right) $ by the elements
in $S_{v}^{sum}$.

 We have then described small subsets of $S$ which suffice to
compute $C^{\geq m}(a)$ and $C^{sum}(a)$. But we still need to show how to
compute these subsets. This is what we do in the next
section.

\section{Algorithms for the conjugacy problem}

 We shall explain in this section our algorithms to compute
 $C^{\geq m}(a)$ and $C^{sum}(a)$, given $a\in G$.

 Let us first explain a technical algorithm, which we did not find in the
 literature. Let $s\in S$ and $v\in M$. We will show how to compute
 their lcm $s\vee v$. More precisely, our algorithm will compute a simple element
 $s'$ such that $s\vee v= v s'$. We must indicate that it is well known how to compute
 the lcm and the gcd of two simple elements, as well as the normal forms of
 any element in $G$.

 \bigskip
 \noindent {\bf \underline{Algorithm 1}} (for computing $s'$ such that $s\vee v= v s'$).
\begin{enumerate}

 \item Compute the normal form of $v= v_1\cdots v_t$.

 \item $s_0=s$.

 \item For every $i=1,\ldots,t$, compute $s_{i-1}\vee v_i$, and write it $v_i s_i$.

 \item Return $s_t$.

\end{enumerate}

\vspace{.3cm}
\begin{proposition}
 Let $s\in S$ and $v\in M$.  Let $s_t$ be the simple element computed by Algorithm~1. Then
 $s\vee v = vs_t$.
\end{proposition}

\begin{proof}
 We proceed by induction on $t=\sup(v)$. If $t=1$ the result is trivial, so suppose that $t>1$
and the result is true for $t-1$. Denote $v'=v_1\cdots v_{t-1}$. We have $s\vee v' =v' s_{t-1}$,
that is, $s_{t-1}$ is the smallest element such that $v's_{t-1}$ is divisible by $s$.
Therefore, an element $r\in M$ satisfies $s\prec vr=v'(v_t r)$ if and only if $s_{t-1}\prec v_t r$,
and this is equivalent to $s_{t-1}\vee v_t  \prec v_t r_t$, that is $v_t s_t\prec v_t r$ hence
$s_t \prec r$. Therefore, $s_t$ is the smallest element satisfying $s\prec v s_t$, as we wanted to show.
\end{proof}

\subsection{Computation of $C^{\geq m}(a)$}

 Let $a\in G$ and $m\in \mathbb{Z}$. As we saw in Section~\ref{secmin},
 the main problem to compute $C^{\geq m}(a)$ is to compute $S_v^{\geq m}$,
 for every $v\in C^{\geq m}(a)$.

 Let then $v\in C^{\geq m}(a)$ and let $x$ be an atom of $M$. Consider the set
 $S_{{\cal P}_v^{\geq m}} \cap mult(x)$. It is always nonempty, since
 $\Delta$ satisfies ${\cal P}_v^{\geq m}$ (for every $v$) and is divisible
 by every atom (by definition of the Garside element).
 Then, by Lemma~\ref{lemrhox}, this set has a unique minimal element, which we denote now
 $r_x$.

 Recall that $S_v^{\geq m}\subset \{r_x: \; x \mbox{ is an atom}\}$, so our first step
 consists of computing $r_x$, for every atom $x\in M$.

  Let $\inf(v)=p\geq m$, so $v=\Delta^m w$ where $w\in M$. Recall that, by
 Proposition~\ref{propcarac}, a simple element $s$ satisfies ${\cal P}_v^{\geq m}$
 if and only if $\tau^m(s)\prec ws$.

 \bigskip
 \noindent {\bf \underline{Algorithm 2}} (for computing $r_x$, minimal element in
 $S_{{\cal P}_v^{\geq m}} \cap mult(x)$).
 \begin{enumerate}

  \item Compute the left normal form of $v=\Delta^p w_1\cdots w_t$.

  \item If $p>m$ then return $x$ and stop.

  \item $w= w_1\cdots w_t$; $\; s=x$.

  \item Compute $\tau^m(s)$.

  \item Use Algorithm~1 to compute $s'$ such that $\tau^m(s) \vee w s = w s s'$.

  \item If $s'=1$ then return $s$ and stop.

  \item $s=ss'$; go to Step 4.

\end{enumerate}

\vspace{.3cm}
\begin{proposition}\label{proalg1}
 Algorithm~2 gives an output, and it is $r_x$.
\end{proposition}

\begin{proof}
 First notice that, if we conjugate any $g \in G$ by a simple element, we can decrease the
 infimum of $g$ by at most one. Hence, if $p>m$, $\inf(x^{-1}vx)\geq m$, so $r_x=x$ and the
 algorithm gives the correct output.

 Now suppose $p=m$, we have computed $w\in M$ such that $v=\Delta^m w$, and we need to
find the smallest $r_x$, such that $x\prec r_x$ and $\tau^m(r_x)\prec wr_x$.
This is done as follows: we take a simple element $s$ such that $x\prec s \prec r_x$
(at the first step $s=x$). Then we use Algorithm 1 to compute $s'$ such that
$\tau^m(s) \vee w s = w s s'$. If $s'=1$ then $\tau^m(s)\prec ws$, so $s=r_x$ and we
obtain the correct output. Otherwise, notice that $\tau^m(s)\prec \tau^m(r_x) \prec w r_x$,
so $w r_x$ is divisible by $\tau^m(s)$ and by $ws$. Hence $wss'\prec wr_x$, so $ss' \prec r_x$.

Therefore, if $s$ is not equal to $r_x$, the algorithm gives an element $s'\neq 1$ such that
$s\prec ss' \prec r_x$, and it starts again checking if $ss'=r_x$. This process must stop,
since the number of left divisors of $r_x$ is finite, so the algorithm finds $r_x$ in finite time.
(moreover if $M$ has homogeneus relations, like $B_n^+$ or $BKL_n^+$, we find $r_x$ in at most
$|\Delta|$ steps).
\end{proof}

\bigskip
\noindent {\bf \underline{Algorithm 3}} (for computing $S_v^{\geq m}$).
\begin{enumerate}

 \item List the atoms of $M$, say $x_1,\ldots, x_{\nu}$. Set $R=\phi$.

 \item For $i=1,\ldots, {\nu}$, do the following:

 \begin{enumerate}

  \item[2a.] Compute $r_{x_i}$, using Algorithm~2.

  \item[2b.] \begin{tabbing} Compute \= $J_i=\{j:\; j\in R \mbox{ and } x_j\prec
  r_{x_i}\}$ and \\  \> $K_i=\{ j: \; j>i \mbox{ and } x_j\prec r_{x_i}\}$.
   \end{tabbing}

  \item[2c.] If $J_i=K_i=\phi$, then set $R=R\cup\{i\}$.
 \end{enumerate}

 \item Return $\{r_{x_i}: \; i\in R\}$.

\end{enumerate}

\vspace{.3cm}
\begin{proposition}
 Algorithm~3 computes $S_v^{\geq m}$.
\end{proposition}

\begin{proof}
 We know by Corollary~\ref{coratoms} that $S_v^{\geq m}$ is the set of minimal elements
in $\{r_{x_i}: \; i=1,\ldots,\nu \}$. We want to find a set $R\subset \{1,\ldots,\nu\}$
such that $S_v^{\geq m}=\{r_{x_i}: \; i\in R \}$. Since we could have $r_{x_i}=r_{x_j}$
for some $i\neq j$, and we want $R$ to be as small as possible, we define it in the following
way: $i\in R$ if and only if $r_{x_i}$ is minimal and there is no $j>i$ such that
$r_{x_i}=r_{x_j}$.


 Suppose that, for some $i$, we already computed the elements in $\{1,\ldots, i-1\}\cap R$
(for $i=1$, this is the empty set). Then we compute $r_{x_i}$ and the sets $J_i$ and $K_i$ defined
in the algorithm. If $i\notin R$, we have two possibilities: either there is some
$j<i$ such that $r_{x_j}$ is a proper divisor of $r_{x_i}$, or there is some $j>i$ such that
$r_{x_j}\prec r_{x_i}$. In the first case $J_i\neq \phi$, and in the latter $K_i\neq \phi$.
Therefore, if both sets are empty, $i\in R$.

 Using this procedure for $i=1,\ldots,\nu$, the algorithm computes $R$, thus $S_v^{\geq m}$.
\end{proof}

\begin{remark}
  Since $S_v^{\geq m}$ is just the subset of $\{r_{x_i}:\;
  i=1,\ldots,\nu\}$ formed by its minimal elements, we could have
  computed $S_v^{\geq m}$ just by comparing the $r_{x_i}$'s and keeping
  the minimal ones. We prefer to use Algorithm~3 since it is
  much faster to see if an atom divides an element, than to
  compare two elements, even if these two elements are simple
  ones.
\end{remark}

\begin{remark}
 Algorithm~3 can still be improved in two different ways. First,
 we do not need to compute all $r_{x_i}$: if during the
 computation of $r_{x_i}$ (using Algorithm~2), we see that
 $x_j\prec s$, for some $j\in R$ or some $j>i$, we can stop
 Algorithm~2 and increase the index $i$ in Algorithm~3. Also, we
 do not need to compute the whole sets $J_i$ and $K_i$: if we find some
 element belonging to one of them, we can directly increase the index $i$.
 We presented Algorithm~3 as above for the clarity of the
 exposition, and because these two improvements do not really
 change the complexity.
\end{remark}

 Finally, let $a\in G$ and $m\in \mathbb{Z}$. The following algorithm works after
 Proposition~\ref{cormainA1}.

 \bigskip
 \noindent {\bf \underline{Algorithm 4}} (for computing $C^{\geq m}(a)$)
\begin{enumerate}

  \item Compute the left normal form of $a$.

  \item Apply repeated cycling to $a$, to obtain $\widehat{a}\in C^{\geq m}(a)$
  (if it exists).

  \item If $\widehat{a}$ is not obtained, return $\phi$ and stop.

  \item Set $v=\widehat{a}$, $\: V=\{\widehat{a}\}$ and $W=\phi$.

  \item Compute $S_v^{\geq m}$, using Algorithm~2.

  \item For every $r\in S_v^{\geq m}$, do the following:

  \begin{enumerate}

     \item[6a.] Set $w=r^{-1}vr\in C^{\geq m}(a)$.

     \item[6b.] Compute the left normal form of $w$.

     \item[6c.] If $w\notin V$, set $V=V\cup \{w\}$.

  \end{enumerate}

  \item $W=W\cup \{v\}$.

  \item If $V=W$ then return $V$ and stop.

  \item Take a new $v\in V\backslash W$; go to Step~5.

\end{enumerate}

\vspace{.3cm}
\begin{remark}
 If we take $m=0$, then $C^{\geq 0}(a)=C^+(a)=C(a)\cap M$, so this algorithm can be used to compute
all positive elements conjugated to $a$.
\end{remark}

\subsection{Computation of $C^{sum}(a)$}

Let now $a\in M$ and $v\in C^{sum}(a)$. Let $\inf(v)=p$ and
$\|v\|=t$  (so $\sup(v)=p+t$), and consider an atom $x$ of $M$.
Similarly to the previous case, the key point to compute
$C^{sum}(a)$ consists on finding the minimal element
$\rho_x\in S_{{\cal P}_v^{sum}}\cap mult(x)$ (it exists since $\Delta$
belongs to this set, and it is unique by Lemma~\ref{lemrhox}).


\bigskip
\noindent {\bf \underline{Algorithm 5}} (for computing $\rho_x$, minimal in
$S_{{\cal P}_v^{sum}}\cap mult(x)$)

 \begin{enumerate}

    \item Using Algorithm~2 compute $r_x$, minimal in $S_{{\cal P}_v^{\geq p}}\cap mult(x)$.

    \item $s=r_x$.

    \item Compute the right normal form of $s^{-1}vs$, say $w_1\cdots w_k\Delta^p$.

    \item If $k=t$, return $s$ and stop.

    \item $s=sw_1$; go to step~3.


 \end{enumerate}

\bigskip
\begin{remark}
 Algorithm~5 can be explained in a very natural way: first compute $r_x$ and the right normal
form $r_x^{-1} v r_x= w_1\cdots w_k \Delta^p$. If $k$ is not minimal, that is, if $k=t+1$, start
{\bf decycling} this right normal form, that is, compute the right normal form of
$w_2\cdots w_k \Delta^p w_1$ (and
multiply $r_x$ by $w_1$). If this new word has not minimal canonical length, decycle again, and
so on. The proposition below shows that this works.
\end{remark}

\begin{proposition}
 Algorithm~5 computes $\rho_x$, the minimal element in
 $S_{{\cal P}_v^{sum}}\cap mult(x)$.
\end{proposition}

\begin{proof}
 We can assume, by Proposition~\ref{proalg1}, that we already know
$r_x$, the minimal element in $S_{{\cal P}_v^{\geq p}}\cap mult(x)$.
We also know that $r_x\prec \rho_x$ by
minimality of $r_x$, since $\rho_x^{-1}v\rho_x\in C^{sum}(a)
\subset C^{\geq p}(a)$. So we have an element $s\in S$ such
that $r_x\prec s \prec \rho_x$, $s^{-1}bs\in C^{\geq p}(a)$
and $\| s^{-1}vs\|\leq t+1$ (at the first step, $s=r_x$).
We compute the right normal form of $s^{-1}vs$,
say $w_1\cdots w_k\Delta^p$.

 There are two possible cases: either $k=t$ or $k=t+1$.
If $k=t$, then $\|s^{-1}vs\|=t$, so $s=\rho_x$ and Step~4 gives the correct
output. Otherwise, $s\neq \rho_x$ and there exists a non trivial element
$s'\in S$ such that $ss'=\rho_x$.
In this case
$$
\rho_x^{-1}v\rho_x= (s')^{-1}s^{-1}vss'=(s')^{-1} w_1\cdots w_{t+1}\Delta^p s'=
u_1\cdots u_t\Delta^p,
$$
where $u_1\cdots u_t\Delta^p$ is in right normal
form. Hence $w_1\cdots w_{t+1}\Delta^p s'= s' u_1\cdots u_t\Delta^p$,
so $w_1\cdots w_{t+1}\tau^{-p}(s')= s' u_1\cdots u_t$.
But $w_1\cdots w_{t+1}$ is in right normal form, so by Corollary~\ref{mainSSS},
$w_1\prec s'$. Therefore $r_x\prec sw_1 \prec \rho_x$, and $s$ is a proper
divisor of $sw_1$. Also, $(sw_1)^{-1}v(sw_1)=w_2\cdots w_{t+1}\Delta^q w_1 \in
C^{\geq p}(a)$ and $\|(sw_1)^{-1}v(sw_1)\|\leq t+1$. So we can set $s=sw_1$, and start again.

 This procedure must stop, finding $s=\rho_x$, since the number of left divisors
of $\rho_x$ is finite. (In homogeneous monoids, the number of steps is bounded by
$|\Delta|$).
\end{proof}

\bigskip
\noindent {\bf \underline{Algorithm 6}} (for computing
$S_v^{sum}$).
\begin{enumerate}

 \item List the atoms of $M$, say $x_1,\ldots, x_{\nu}$. Set $R=\phi$.

 \item For $i=1,\ldots, \nu$, do the following:

 \begin{enumerate}

  \item[2a.] Compute $\rho_{x_i}$, using Algorithm~5.

  \item[2b.] \begin{tabbing} Compute \= $J_i=\{j:\; j\in R \mbox{ and } x_j\prec
  \rho_{x_i}\}$ and  \\  \> $K_i=\{ j: \; j>i \mbox{ and } x_j\prec
  \rho_{x_i}\}$. \end{tabbing}

  \item[2c.] If $J_i=K_i=\phi$ then $R=R\cup\{i\}$.
 \end{enumerate}

 \item Return $\{\rho_{x_i}: \; i\in R\}$.

\end{enumerate}

\begin{proposition}
 Algorithm~6 computes $S_v^{sum}$.
\end{proposition}

\begin{proof}
This algorithm parallels Algorithm~3, and it
works in the same way, since $S_v^{sum}$ is the set of minimal
elements in $\{\rho_{x_i}:\; i=1,\ldots,n \}$.
\end{proof}

Finally, the following algorithm is analogous to Algorithm~4.
Using the previous algorithms, it will then compute $C^{sum}(a)$ for
any given $a\in G$.

 \bigskip
 \noindent {\bf \underline{Algorithm 7}} (for computing $C^{sum}(a)$)
\begin{enumerate}

  \item Compute the left normal form of $a$.

  \item Using cyclings and decyclings, compute $\widetilde{a}\in
           C^{sum}(a)$.

  \item $v=\widetilde{a}$; $\: V=\{\widetilde{a}\}$;\: $W=\phi$.

  \item Compute $S_v^{sum}$, using Algorithm~6.

  \item For every $r\in S_v^{sum}$, do the following:

  \begin{enumerate}

     \item[5a.] $w=r^{-1}vr\in M$.

     \item[5b.] Compute the left normal form of $w$.

     \item[5c.] If $w\notin V$ then $V=V\cup \{w\}$.

  \end{enumerate}

  \item $W=W\cup \{v\}$.

  \item If $V=W$ then return $V$ and stop.

  \item Take a new $v\in V\backslash W$; go to Step~4.

\end{enumerate}

\section{Complexity}

  In this section we shall study the complexity of our algorithms,
applied to several examples of Garside monoids, such as $B_n^+$,
$BKL_n^+$ and Artin monoids. We do not discuss here the general
case, since the complexity strongly depends on the way of computing
normal forms, LM($a$), $a\vee b$, $\widetilde{a}$, etc. in each particular
case.

 We shall give theoretical upper bounds for this complexity. In the next section we will
also compare, with many examples in $B_n^+$, the running time of our algorithm with the one
by Elrifai and Morton, to show that our improvement is significant in practice.

 We know that the theoretical results in this section can be improved: we are just
interested on showing that our algorithms are much faster than the preceding ones, while
sharper bounds for the complexity would require a deeper study.

\subsection{The Artin braid monoid $\mathbf{B_n^+}$.}

\noindent
{\bf \underline{Complexity of computing $\mathbf{C^{\geq m}(a)}$}}

\vspace{.2cm}
Recall the definition of $B_n^+$, and notice that the relations are
homogeneous, so every two conjugate elements in $B_n^+$ have the same word
length. Let then $v\in B_n^+$ be of word length $l$. Its left normal form
$v_1\cdots v_t$ verifies $t\leq l$.

\vspace{.3cm}
\noindent {\bf Algorithm~1.} It computes the normal form of $v$, which takes time
$O(l^2n\log n)$ (see\cite{T}). Then it computes $t$ words, $s_1,\ldots, s_t$.
The computation of each one takes time $O(n\log n)$ (\cite{T}), so all of them
are computed in time $O(l n\log n)$, and the whole algorithm has complexity
$O(l^2n\log n)$.

 Now suppose that $v\in B_n$, and its left normal form is $\Delta^p w_1\cdots w_t$.
 One has $t\leq l$. In $B_n$, the homomorphism $\tau$ can be defined as follows:
$\tau(\sigma_i)=\sigma_{n-i}$. Hence $\tau^2=\mbox{id}$, so for every word $w$,
$\tau^m(w)$ can be computed in time $O(|w|)$.

\vspace{.3cm}
\noindent {\bf Algorithm~2.} First it computes the left normal form of $v$ in time
 $O(l^2 n\log n)$. If $p=m$, it runs the loop consisting of Steps 4-7.
 The two important steps are the following:
\begin{itemize}
 \item {\bf Step 4:} Notice that, every time the loop is repeated, we already know
 $\tau^m(s)$ for the old value of $s$. So in order to compute
 $\tau^m(ss')=\tau^m(s)\tau^m(s')$, we just need to compute $\tau^m(s')$, which is $O(|s'|)$.
 Since the product of all possible $s'$ is still a simple element, $r_x$, all repetitions of this step
 can be made in  time $O(n(n-1)/2)$, which is the length of $\Delta$.

 \item {\bf Step 5:}
 Here, when we apply Algorithm~1, we compute the normal form of $ws$, and the elements
 $s_1,\ldots,s_t$, where $s$ runs over an ascending chain of divisors of $\Delta$.
 As above, all these computations together require the same number of operations as
 just applying Algorithm~1 to $wr_x$ (see~\cite{T}). Moreover, we have computed the normal
 form of $w$ at the beginning of Algorithm~2 so, again by~\cite{T}, the normal form of $wr_x$
 can be computed in time $O(ln\log n)$. Hence all repetitions of this step can be
 made in time $O(ln\log n)$.
\end{itemize}

 Therefore, the theoretical complexity of Algorithm~2 is $O(l^2n^2)$.

\vspace{.3cm}
\noindent {\bf Algorithm~3.} The only non-negligible step is the
loop of Step~2, which is repeated $n-1$ times (the number of atoms in $B_n^+$),
and does the following: it computes $r_{x_i}$ (which takes time
$O(l^2n^2)$) and it verifies at most $n-2$ times if an atom divides
$r_{x_i}$ (this takes time $O(n)$ by \cite{T}). Hence, the complexity of
Algorithm~3 is $O(l^2n^3)$.

\vspace{.3cm}
 Finally, let $a\in B_n$ be given as a word of length $l$ in the generators $\sigma_i$.

\vspace{.3cm}
\noindent {\bf Algorithm~4.} It starts by computing the normal form
of $a$ (time $O(l^2n\log n)$). Then it finds $\widehat{a}$, which takes $O(l^2n^3)$
by~\cite{BKL2}. Next it starts a loop, which is repeated $k$ times
(the number of elements in $C^{\geq m}(a)$), and does the following: First, it computes
$S_v^{\geq m}$, taking time $O(l^2n^3)$. Then it runs another loop,
repeated at most $n-1$ times, which works at follows:
\begin{itemize}

\item It computes the normal form of a word of length $l$, thus taking time
  $O(l^2n\log n)$.

\item It verifies if an element is in a list $V$, taking a negligible time compared
to the previous step.
\end{itemize}

Finally, it verifies if $V=W$, but since $W\subset V$ we just have to compare
the lengths. The time to do this is negligible. Therefore, the complexity of
Algorithm~4 (i.e. the complexity of computing $C^{\geq m}(a))$, is
$O\left(l^2n\log n + l^2 n^3+k(l^2n^3+(n-1)l^2n\log n)\right)$,
which yields the following.

\begin{proposition}
 Given $a\in B_n$ as a word of length $l$, the complexity of computing $C^{\geq m}(a)$
 (for the Artin presentation) is $O\left(kl^2n^3\right)$, where $k$ is the number of elements
 in $C^{\geq m}(a)$.
\end{proposition}

 Remark that if we try to compute $C^{\geq m}(a)$ using the techniques of Elrifai and Morton,
we should use Algorithm~4, but replacing $S_v^{\geq m}$ by the whole $S$, which has cardinality
$n!$. The time would be in this case $O\left(kl^2(n!)n\log n\right)$.

\bigskip
\noindent
{\bf \underline{Complexity of computing $\mathbf{C^{sum}(a)}$}}

\vspace{.2cm}
The study of the complexity of Algorithms~5, 6 and 7, is very similar to that of
Algorithms~2, 3 and 4.

\vspace{.3cm}
\noindent {\bf Algorithm~5.} It starts by computing $r_x$, taking $O(l^2n^2)$.
Next it computes a right normal form ($O(l^2n\log n)$), and then it does a number of decyclings,
which is bounded by $\frac{n(n-1)}{2}$. By~\cite{BKL2}, each decycling takes time $O(ln)$,
so the whole complexity of Algorithm~5 is $O(l^2n^3)$.

\vspace{.3cm}
\noindent {\bf Algorithm~6.} It does the same as Algorithm~3, but it computes $\rho_{x_i}$
instead of $r_{x_i}$ in Step~2a. Hence its complexity is $O(l^2n^4)$.

\vspace{.3cm}
\noindent {\bf Algorithm~7.} It has two main differences with respect to Algorithm~4.
It computes $\widetilde{a}$ instead of $\widehat{a}$ (but this can be made in $O(l^2n^3)$
by~\cite{BKL2}), and it computes $S^{sum}_v$ instead of $S^{\geq m}_v$. Therefore one has the
following.

\begin{proposition}
 Given $a\in B_n$ as a word of length $l$, the complexity of computing $C^{sum}(a)$ (for
 the Artin presentation) is
$O(kl^2n^4)$, where $k$ is the number of elements in $C^{sum}(a)$.
\end{proposition}

Remark that, if we compute the complexity of the algorithm by Elrifai and Morton, using the
above methods, we obtain $O\left(kl^2(n!)n\log n\right)$.

\subsection{The Birman-Ko-Lee monoid $\mathbf{BKL_n^+}$.}

\noindent
{\bf \underline{Complexity of computing $\mathbf{C^+(a)}$}}

\vspace{.2cm}
 We just need to follow here the same reasoning that in the previous subsection, taking into
account the differences between $BKL_n^+$ and $B_n^+$. The complexities of the basic computations
in $BKL_n^+$ can be found in \cite{B-K-L}. For instance, in $BKL_n^+$, the computation of each
$s_i$ in Algorithm~1 takes time $O(n)$, the length of the Garside element
is $n-1$, and the normal form of a word $w$ is computed in time $O(|w|^2 n)$.
This implies that Algorithm~1 and Algorithm~2 both have complexity $O(l^2n)$.

 In order to study Algorithm~3, we must know that the number of atoms in $BKL_n^+$ is
$\frac{n(n-1)}{2}$, and to check if an atom divides a simple element takes time $O(n)$,
so Algorithm~3 has complexity $O(l^2n^5)$.

 Finally, the computation of $\widehat{a}$ takes time $O(l^2n^2)$ (see \cite{BKL2}),
so the above method to compute the complexity of Algorithm~4 yields the following:

\begin{proposition}
 Given $a\in BKL_n$ as a word of length $l$, the complexity of computing $C^{\geq m}(a)$
 (for the Birman-Ko-Lee presentation) is
$O\left(kl^2n^5 \right)$, where $k$ is the number of elements in $C^{\geq m}(a)$.
\end{proposition}

 As above, the complexity of the previously known algorithm is much worse:
$O(kl^2{\cal C}_nn)$, where $3^n<{\cal C}_n<4^n$.

\bigskip
\noindent
{\bf \underline{Complexity of computing $\mathbf{C^{sum}(a)}$}}

\vspace{.2cm}
 We just need to know that the computation of $\widetilde{a}$ takes time $O(l^2n^2)$ (see \cite{BKL2}).
 Hence the complexities of Algorithms~5, 6 and 7 are respectively $O(l^2n^2)$, $O(l^2n^5)$ and
 $O(kl^2n^5)$. Therefore, one has:

\begin{proposition}
Given $a\in BKL_n$ as a word of length $l$, the complexity of computing $C^{sum}(a)$
(for the Birman-Ko-Lee presentation) is $O(kl^2n^5)$, where $k$ is the number of elements
in $C^{sum}(a)$.
\end{proposition}

 Notice that the complexity of the known algorithm was $O(kl^2{\cal C}_nn)$,
so our algorithm improves it considerably.

 One interesting remark is that our algorithm works faster, a priori, for the monoid $B_n^+$ than for
$BKL_n^+$. This is due to a simple fact: in our algorithm the number of atoms is more relevant than
the number of simple elements. In $BKL_n^+$, the number of simple elements is much smaller than in
$B_n^+$, but the number of atoms is $\frac{n(n-1)}{2}$, while in $B_n^+$ is $n-1$.

\subsection{Artin monoids}

 As we mentioned in the introduction, the Artin groups of finite type are Garside groups, so we can apply
our algorithms to the corresponding Artin monoids (see~\cite{B}  for an introduction to Artin monoids
and groups). In \cite{BS} we can find algorithms to deal with Artin monoids: computation of normal forms,
greatest common divisors, division algorithms, etc. Although these algorithms seem to be exponential in the
length of the words involved, in~\cite{C} it is shown that finite type Artin groups are biautomatic, so
there are quadratic algorithms to compute all of the above.

 Nevertheless, since we are mainly interested in comparing our algorithms with the previous ones, we just need
to know the length of the Garside element $\Delta$, and the number of simple elements in any given Artin group.
Let then $G$ be an Artin group of rank $n$, that is, $A_n$, $B_n$, $D_n$, $E_n$ (if $n=6,7,8$), $F_n$ (if $n=4$),
$H_n$ ($n=3,4$) or $I_2(p)$ (if $n=2$), and let $h$ be its Coxeter number. It is known that $|\Delta|=\frac{nh}{2}$,
where $h=O(n)$, and that $\#(S)\geq n!$.

  Hence, if the complexity of the conjugacy algorithm by Elrifai and Morton is $O(xn!)$ for some $x$
depending on $n$ and $l$, our algorithm will have complexity $O(xn^3)$. This is shown by using the same
arguments as in the previous subsections.

\section{Effective computations}

\subsection{Comparison with the Elrifai-Morton algorthim}

 In the previous section, we found theoretical upper bounds for the complexity of our
algorithms. We showed that our algorithm is, in theory, much better than the Elrifai-Morton
one (for $n>5$). In this section we effectively compare the two algorithms, in the following way:
For given $n$ and $l$, ($3\leq n \leq 5$ and $10\leq l \leq 20$) we took 5000 random pairs of
positive braids in $B_n$ of length $l$ (using Artin presentation), we tested conjugacy using
both algorithms, and we compared the Average Running Time (ART) and the Maximum Running Time (MRT).
We did the same for $n=6$ and $l=10$, for 1144 pairs.

 We can conclude that our algorithm is faster for $n\geq 4$, and much faster for $n\geq 5$ (We were not
able to compute the cases $n=5$ and $l=19,20$ using the Elrifai-Morton algorithm since the computations
were too long).

In the tables below one can see the results: We wrote F-GM for our algorithm and E-M for the
Elrifai-Morton one. The time is given in seconds.

 {\small \vspace{.2cm} \noindent
$$n=3$$
\noindent
\begin{tabular}{|c||c|c|c|c|c|c|}
\hline $l$ & 10 & 11 & 12 & 13 & 14 & 15    \\
\hline \hline ART F-GM & 0.1526 & 0.2011 &  0.2361 & 0.3038 &
0.3386 & 0.3951   \\
\hline ART E-M & 0.1144 & 0.1460 & 0.1692 & 0.2133 & 0.2367 &
0.2723  \\ \hline
\hline MRT F-GM & 2.429 & 3.599 & 4.680 & 6.080 & 7.450 & 6.960  \\
\hline MRT E-M & 1.659 & 2.539 & 3.220 & 4.089 & 5.029 & 4.599
  \\
\hline
\end{tabular}

\vspace{.5cm}
\hfill
\begin{tabular}{|c||c|c|c|c|c|}
\hline $l$ &  16 & 17 & 18 & 19
& 20    \\
\hline \hline ART F-GM &  0.3896 &  0.5021 & 0.5473 & 0.6494 &
0.7292   \\
\hline ART E-M &  0.2710 & 0.3392 & 0.3710 & 0.4329 & 0.4841 \\
\hline   \hline MRT F-GM &
11.299 & 10.530 &
12.469 & 15.090 & 16.539 \\
\hline MRT E-M &
7.219 & 6.729 & 7.970 &  9.950 & 11.039
  \\
\hline
\end{tabular}

\vspace{0.8cm}

\noindent
$$n=4$$
\noindent
\begin{tabular}{|c||c|c|c|c|c|c|}

\hline $l$ & 10 & 11 & 12 & 13 & 14 & 15     \\
\hline \hline ART F-GM & 0.3559 & 0.4796 & 0.6772 & 0.7870 & 1.0264
& 1.2599   \\
\hline ART E-M &      0.6118 &
               0.7233 &
               1.0127 &
               1.2086 &
               1.5909 &
               1.9538   \\ \hline
\hline MRT F-GM & 8.680 & 11.390 & 16.519 & 23.949 & 33.969 & 42.029  \\
\hline MRT E-M &   16.319 &
                  22.329 &
                  28.440 &
                  41.579 &
                  61.790 &
                  74.999
  \\
\hline
\end{tabular}

\vspace{.5cm}
\hfill
\begin{tabular}{|c||c|c|c|c|c|}

\hline $l$ &  16 & 17 & 18 & 19
& 20    \\
\hline \hline ART F-GM &  1.4548 &
1.7436 & 2.2029 & 2.6616 & 2.9942   \\
\hline ART E-M &
               2.3106 &
               2.7995 &
               3.5548 &
               4.3280 &
               4.7226   \\ \hline
\hline MRT F-GM &  41.910 &
62.940 & 72.940 & 103.470 & 148.989 \\
\hline MRT E-M &
                  70.039 &
                 107.969 &
                 137.720 &
                 173.060 &
                 245.740
  \\
\hline
\end{tabular}

\vspace{.8cm}

\noindent
$$n=5$$
\noindent
\begin{tabular}{|c||c|c|c|c|c|c|}

\hline $l$ & 10 & 11 & 12 & 13 & 14 & 15    \\
\hline \hline ART F-GM & 1.0997 &
               1.8463 &
               2.7657 &
               3.7962 &
               3.8195 &
               4.4797  \\
\hline ART E-M & 7.8690 & 11.1207 & 17.1455 & 23.2491 & 26.2595 & 29.7934    \\
\hline \hline MRT F-GM &
 21.239 &
                  46.070 &
                  65.530 &
                  88.940 &
                  139.180 &
                 155.260
 \\
\hline MRT E-M &177.489&322.039&456.669&611.609 & 1068.970 & 1178.790  \\
\hline
\end{tabular}

\vspace{.5cm}
\hfill
\begin{tabular}{|c||c|c|c|c|c|}
\hline $l$ &  16 & 17 & 18 & 19
& 20    \\
\hline \hline ART F-GM &
               5.6410 &
               7.1540 & 8.8198 & 9.4597 & 10.6614  \\
\hline ART E-M &  38.7974 & 51.0028 & 62.0018 & &    \\
\hline \hline MRT F-GM &

                 254.770 &
                 411.320 & 401.409 & 516.119& 532.469
 \\
\hline MRT E-M & 2116.239 & 3221.880 & 3218.93 & &   \\
\hline
\end{tabular}

\vspace{.8cm}

\quad \quad
$n=6$

\medskip
\noindent
\begin{tabular}{|c||c|}

\hline $l$ & 10  \\
\hline \hline ART F-GM & 2.2450    \\
\hline ART E-M & 506.224     \\
\hline \hline MRT F-GM & 43.935
 \\
\hline MRT E-M & 7495.288   \\
\hline
\end{tabular}

 }

\subsection{Exhaustive computation of conjugacy classes and summit classes}

 In the previous section, we saw that the complexity of all our algorithms depends on the size of the sets
$C^{\geq m}(a)$ or $C^{sum}(a)$, for $a\in M$. In the cases of $B_n^+$ or $BKL_n^+$, the only upper bounds known
for these sets are exponential in $n$ and in $l=|a|$. Nevertheless, we have the following (recall that,
in this case, $C^+(a)=C^{\geq 0}(a)=C(a)\cap B_n^+$):

\bigskip
\noindent
{\bf Conjecture:} (Thurston, \cite{T}) Let $n$ be a fixed integer and let $a\in B_n^+$, having word length $l$.
There is an upper bound for $C^+(a)$ which is a polynomial in $l$.

\bigskip
The existence of this upper bound for $C^+(a)$, or even for $C^{sum}(a)$, would imply the following:

\bigskip
\noindent
{\bf Conjecture:} (Birman, Ko and Lee, \cite{BKL2}) For every fixed integer $n$, there exists a solution for
the conjugacy problem in $B_n$, which is polynomial in the word length of the elements involved.

\bigskip
 In order to have some numerical evidence to support these conjectures, we have computed,
for $n=3,\ldots,8$ and several values of $l$, all the conjugacy classes
of words of length $l$ in $B_n^+$, as well as the corresponding summit classes.
In the tables below we present the following data, for the set $W_l$ of elements in $B_n^+$ having word
length $l$:
\begin{itemize}

\item $\mathbf{CC^+}$: The number of Conjugacy Classes in $W_l\subset B_n^+$.

\item {\bf max }$\mathbf{C^+}$: The size of the biggest one. That is, the number of elements
in the biggest $C^+(a)$, for $a\in W_l$.

\item {\bf max }$\mathbf{C^{sum}}$: The size of the biggest summit class.

\item $\mathbf{v}$:  A representative from one of those biggest summit class. That is, an element
     $v\in C^{sum}(a)$, where $C^{sum}(a)$ has maximal size.

\end{itemize}

%

\bigskip
\centerline{
\begin{tabular}{|c|c|c|c|c|}
\hline
\multicolumn{5}{|c|}{\bf n=3} \\
 \hline
 $l$  & $CC^+$ & $\max C^+$ &  $\max C^{sum}$ & $v$ \\
\hline
\hline
 4   & 3     & 6   & 2  &  $\sigma_1^3\sigma_2 $ \\
\hline
 5   & 3     & 10   & 6   &  $\sigma_1^3\sigma_2^2 $\\
\hline
  6   & 5     & 12   & 8  &  $\sigma_1^4\sigma_2^2 $ \\
\hline
  7   & 5     & 16   & 10  & $\sigma_1^5\sigma_2^2$ \\
\hline
  8   & 8     & 20   & 12  & $\sigma_1^6\sigma_2^2$ \\
\hline
  9   & 9     & 29   & 14  & $\sigma_1^7\sigma_2^2$ \\
\hline
  10  & 13    & 30   & 16  & $\sigma_1^8\sigma_2^2$ \\
\hline
   11  & 16    & 40   & 18  & $\sigma_1^{9}\sigma_2^2$ \\
\hline
   12   & 27     & 48   & 20  & $\sigma_1^{10}\sigma_2^2$ \\
\hline
   13   & 33     & 64   & 22  & $\sigma_1^{11}\sigma_2^2$ \\
\hline
   14   & 50     & 80   & 24  & $\sigma_1^{12}\sigma_2^2$ \\
\hline
   15   & 70     & 125   & 26  & $\sigma_1^{13}\sigma_2^2$ \\
\hline
   16   & 107     & 126   & 28  & $\sigma_1^{14}\sigma_2^2$ \\
\hline
   17   & 153     & 160   & 30  & $\sigma_1^{15}\sigma_2^2$ \\
\hline
   18  & 241     & 192   & 32  & $\sigma_1^{16}\sigma_2^2$ \\
\hline
   19  & 349     & 256   & 34  & $\sigma_1^{17}\sigma_2^2$ \\
\hline
   20  & 542     & 320   & 36  & $\sigma_1^{18}\sigma_2^2$ \\
\hline
\end{tabular}
}

\bigskip
\centerline{
\begin{tabular}{|c|c|c|c|c|}
\hline
\multicolumn{5}{|c|}{\bf n=4} \\
\hline
 $l$  & $CC^+$ & $\max C^+$ &  $\max C^{sum}$ & $v$ \\
\hline
\hline
 4 & 7 & 12 & 4 & $\sigma_1^3 \sigma_2$ \\
\hline
 5 & 9 & 20 & 12 & $\sigma_1^3 \sigma_2^2$ \\
\hline
 6 & 16 & 40 & 16 & $\sigma_1^4 \sigma_2^2$ \\
\hline
 7 & 21 & 54 & 22 & $\sigma_1^5 \sigma_2 \sigma_3$ \\
\hline
 8 & 36 & 72 & 32 & $\sigma_1^6 \sigma_2 \sigma_3$\\
\hline
 9 & 54 & 94 & 50 & $\sigma_1^4 \sigma_2^2 \sigma_3^2 \sigma_2$\\
\hline
 10 & 96 & 156 & 60 & $\sigma_1^5 \sigma_2^2 \sigma_3^2 \sigma_2$\\
\hline
 11 & 160 & 252 & 70 & $\sigma_1^6 \sigma_2^2 \sigma_3^2 \sigma_2$\\
\hline
 12 & 304 & 344 & 88 & $\sigma_1^5\sigma_2^2\sigma_1\sigma_3\sigma_1\sigma_2\sigma_3$\\
\hline
 13 & 538 & 582 & 114 & $\sigma_1^6\sigma_2^2\sigma_1\sigma_3\sigma_1\sigma_2\sigma_3$\\
\hline
 14 & 1030 & 752 & 140 & $\sigma_1^7\sigma_2^2\sigma_1\sigma_3\sigma_1\sigma_2\sigma_3$\\
\hline
 15 & 1954 & 1114 & 166 & $\sigma_1^8\sigma_2^2\sigma_1\sigma_3\sigma_1\sigma_2\sigma_3$\\
\hline
\end{tabular}
}

\bigskip
\centerline{
\begin{tabular}{|c|c|c|c|c|}
\hline
\multicolumn{5}{|c|}{\bf n=5} \\
\hline
 $l$  & $CC^+$ & $\max C^+$ &  $\max C^{sum}$ & $v$ \\
\hline
\hline
 4 & 10 & 24 & 8 & $\sigma_1^2 \sigma_2 \sigma_3$ \\
\hline
 5 & 15 & 36 & 18 & $\sigma_1^3 \sigma_2^2$\\
\hline
 6 & 28 & 80 & 24 & $\sigma_1^4 \sigma_2^2$ \\
\hline
 7 & 44 & 136 & 44 & $\sigma_1^5 \sigma_2 \sigma_3$ \\
\hline
 8 & 81 & 188 & 64 & $\sigma_1^6 \sigma_2 \sigma_3$ \\
\hline
 9 & 141 & 288 & 104 & $\sigma_1^5 \sigma_2 \sigma_3 \sigma_2 \sigma_4$ \\
\hline
 10 & 281 & 516 & 156 & $\sigma_1^6 \sigma_2 \sigma_3 \sigma_2 \sigma_4$ \\
\hline
 11 & 520 & 702 & 208 & $\sigma_1^7 \sigma_2 \sigma_3 \sigma_4^2$ \\
\hline
 12 & 1194 & 1018 & 260 & $\sigma_1^8 \sigma_2 \sigma_3 \sigma_4^2$\\
\hline
\end{tabular}
}

\bigskip
\centerline{
\begin{tabular}{|c|c|c|c|c|}
\hline
\multicolumn{5}{|c|}{\bf n=6} \\
\hline
 $l$  & $CC^+$ & $\max C^+$ &  $\max C^{sum}$ & $v$ \\
\hline
\hline
 4 & 13 & 36 & 16 & $\sigma_1 \sigma_2 \sigma_3 \sigma_4$\\
\hline
 5 & 22 & 56 & 30 & $\sigma_1 \sigma_2 \sigma_1 \sigma_4^2$\\
\hline
 6 & 44 & 120 & 36 & $\sigma_1^4 \sigma_2 \sigma_3 $  \\
\hline
 7 & 76 & 272 & 72 & $\sigma_1^4 \sigma_2 \sigma_3 \sigma_4$\\
\hline
 8 & 148 & 412 & 124 & $\sigma_1^5 \sigma_2 \sigma_3 \sigma_4$\\
\hline
 9 & 276 & 576 & 208 & $\sigma_1^5 \sigma_2 \sigma_3 \sigma_4^2$\\
\hline
 10 & 573 & 1032 & 372 & $\sigma_1^5 \sigma_2 \sigma_3 \sigma_4 \sigma_5^2$\\
\hline
\end{tabular}
}

\bigskip
\centerline{
\begin{tabular}{|c|c|c|c|c|}
\hline
\multicolumn{5}{|c|}{\bf n=7} \\
\hline
 $l$  & $CC^+$ & $\max C^+$ &  $\max C^{sum}$ & $v$ \\
\hline
\hline
 4 & 14 & 60 & 24 & $\sigma_1 \sigma_2 \sigma_4^2$\\
\hline
 5 & 26 & 84 & 60 & $\sigma_1 \sigma_2 \sigma_1 \sigma_4^2$\\
\hline
 6 & 56 & 160 & 72 & $\sigma_1 \sigma_2 \sigma_4 \sigma_2 \sigma_1^2$\\
\hline
 7 & 104 & 408 & 108 & $\sigma_1^4 \sigma_2 \sigma_3 \sigma_4$\\
\hline
 8 & 215 & 824 & 192 & $\sigma_1^4 \sigma_2 \sigma_3 \sigma_4 \sigma_5$\\
\hline
 9 & 424 & 1160 & 416 & $\sigma_1^4 \sigma_2 \sigma_3 \sigma_4 \sigma_5^2$\\
\hline
 10 & 914 & 1992 & 744 & $\sigma_1^5 \sigma_2 \sigma_3 \sigma_4 \sigma_5^2$\\
\hline
\end{tabular}
}

\bigskip
\centerline{
\begin{tabular}{|c|c|c|c|c|}
\hline
\multicolumn{5}{|c|}{\bf n=8} \\
\hline
 $l$  & $CC^+$ & $\max C^+$ &  $\max C^{sum}$ & $v$ \\
\hline
\hline
 4 & 15 & 100 & 48 & $\sigma_1 \sigma_2 \sigma_3 \sigma_5$ \\
\hline
 5 & 29 & 144 & 100 & $\sigma_1 \sigma_2 \sigma_1 \sigma_4^2$\\
\hline
 6 & 66 & 216 & 144 & $\sigma_1 \sigma_2 \sigma_1 \sigma_3 \sigma_5^2$\\
\hline
 7 & 130 & 544 & 168 & $\sigma_1 \sigma_2 \sigma_1 \sigma_3 \sigma_4 \sigma_6^2$\\
\hline
 8 & 281 & 1236 & 360 & $\sigma_1 \sigma_2 \sigma_1 \sigma_4^3 \sigma_5^2$\\
\hline
\end{tabular}
}

\begin{acknowledge}
 The main ideas in this work were developed during a stay of both authors at the {\em Laboratoire
de Topologie de l'Universit\'e de Bourgogne} at Dijon (France). We are very grateful to all members of
the Laboratoire, and in particular to Luis Paris for his many useful suggestions.
We are also grateful to Jean Michel for his precise and helpful comments on an earlier version of this
paper.
\end{acknowledge}

\end{document}